%% file: xxx_bergen_proceedings.tex
\RequirePackage{snapshot}

\documentclass[11pt, a4paper]{article}

\usepackage{amsmath}
    \allowdisplaybreaks
\usepackage{amssymb}
\usepackage{amsthm}  
    \newtheorem{thm}{Theorem}
    \newtheorem{prop}[thm]{Proposition}

    \theoremstyle{definition}
    \newtheorem{defn}{Definition}

    \theoremstyle{definition}
    
\usepackage[mathic=true]{mathtools} 
\usepackage{etoolbox}

\usepackage{graphicx}
    \newcommand{\figdir}{pstex}
    
\usepackage{graphbox} 

\usepackage{subcaption} 

\usepackage{pgfplots}
    \pgfplotsset{compat=1.18}
    \usetikzlibrary{shapes,arrows}
    \usetikzlibrary{mindmap}
    \usetikzlibrary{decorations.fractals, positioning}

    \tikzstyle{point} = [circle, draw, fill=blue!80, 
        text width=0.5em, text centered, node distance=3cm, inner sep=0pt]
    \tikzstyle{block} = [rectangle, draw, fill=blue!20, text width=5em, text centered, rounded corners, minimum height=2em]

    \tikzstyle{line} = [draw, -latex']

\usepackage{quiver}

\usepackage[american, british]{babel}
\usepackage[english=british]{csquotes}

\usepackage[style=ieee, backref=true, 
    backend=bibtex,
    hyperref=true, giveninits=true, sorting=nyt, language=british]{biblatex}
    \DefineBibliographyStrings{english}{%
        backrefpage = {cited on p.~},
        backrefpages = {cited on pp.~},
        url = {},
    } 
    \DefineBibliographyExtras{british}{%
        \stdpunctuation
    }

\usepackage{xpatch}
    \xpatchfieldformat{eprint:arxiv}{\space}{}{}{}
    \xpatchfieldformat{doi}{\space}{}{}{}
    \xpatchfieldformat{url}{\addcolon\space}{}{}{}

\usepackage[T1]{fontenc}
\usepackage[largesc]{newpxtext}
\usepackage[vvarbb]{newpxmath}

\usepackage{anyfontsize} 

\usepackage{accents}
\usepackage{stmaryrd} 

\usepackage{booktabs}

\usepackage[immediate]{silence} 
\usepackage{sectsty}
    \chapternumberfont{\sffamily\bfseries\LARGE\sectionrule{3ex}{3pt}{-1ex}{1pt}}
    \sectionfont{\rmfamily\bfseries\sectionrule{3ex}{0pt}{-1ex}{1pt}}
    \subsectionfont{\rmfamily\bfseries\large}
    \subsubsectionfont{\rmfamily\bfseries}

    \DeactivateWarningFilters[temp]

\usepackage[inline, shortlabels]{enumitem}

\numberwithin{equation}{subsection}

\usepackage[unicode, naturalnames, colorlinks=true]{hyperref}

\usepackage{microtype}

\newcommand{\catstyle}{\mathcal}
\newcommand{\outercat}{\mathbf}

\newcommand{\catname}[1]{\csdef{#1}{\outercat{\catstyle #1}}}
\catname{Set}
\catname{Vect}
\catname{Ab}
\catname{Rep}
\catname{Top}
\catname{Act}
\catname{Cat}
\catname{Mod}
\catname{V}
\catname{W}
\catname{C}
\catname{D}
\catname{B}
\catname{O}
\newcommand{\FinVect}{\outercat{\catstyle{V}ect_{fd}}}

\newcommand{\selfenr}[1]{\undertilde{#1}}

\newcommand{\Rbar}{\overline{\mathbb{R}}}
\newcommand{\Rplusbar}{\overline{\mathbb{R}}_{+}}

\newcommand{\terminalset}{\{\star\}}

\newcommand{\Q}{\mathbb{Q}}

\newcommand{\Z}{\mathbb{Z}}
\newcommand{\CC}{\mathbb{C}}
\newcommand{\F}{\mathbb{F}}

\newcommand{\define}{\textbf}
\newcommand{\formal}[1]{\ulcorner\kern-4\scriptspace #1\kern-4\scriptspace\urcorner}

\newcommand{\ob}[1]{\operatorname{ob}(#1)}
\newcommand{\op}{\mathrm{op}}
\newcommand{\one}{\mathbb{1}}
\newcommand{\R}{\mathbb{R}}

\DeclareMathOperator{\Fix}{Fix}

\renewcommand{\phi}{\varphi}
\newcommand{\id}{\operatorname{id}}
\newcommand{\dd}{\operatorname{d}}
\newcommand{\minusdot}{\stackrel{.}{-}}

\newcommand{\profto}{\nrightarrow}

\DeclarePairedDelimiterX{\funcat}[2]{\llbracket}{\rrbracket}{#1 , \mathopen{} #2}
\DeclarePairedDelimiterXPP{\cvx}[2]{\operatorname{Cvx}}{(}{)}{}{#1 , \mathopen{} #2}
\DeclarePairedDelimiterX{\closedinterval}[2]{[}{]}{#1 , \mathopen{} #2}
\DeclarePairedDelimiterXPP{\innerprod}[3]{}{\langle}{\rangle}{_{#3}}{#1, \mathopen{} #2}
\DeclarePairedDelimiter{\card}{\lvert}{\rvert}
\DeclarePairedDelimiter{\size}{\lvert}{\rvert}
\DeclarePairedDelimiter{\magnitude}{\lvert}{\rvert}

\newcommand{\hausdorff}{\operatorname{d}_{\mathrm{H}}}
\newcommand{\Fun}{\operatorname{Fun}}
\newcommand{\Hom}{\operatorname{Hom}}

\renewcommand{\coloneq}{\coloneqq}

\DeclareMathOperator{\Cvx}{Cvx}

\DeclareMathOperator{\fenchpro}{\mathbb{L}}
\DeclareMathOperator{\fenchb}{\mathbb{L}^\ast}
\DeclareMathOperator{\rfenchb}{\mathbb{L}_\ast}

\newcommand{\fenchbop}[1]{\fenchb(#1)}
\newcommand{\rfenchbop}[1]{\rfenchb(#1)}

\newcommand{\dual}[1]{{#1}^{\#}}

\title{Metric-like spaces as enriched categories:\\ three vignettes}
\author{Simon Willerton}
\date{}

\begin{document}

\maketitle

\begin{abstract}
    This is a write-up of a talk given at the CATMI meeting in Bergen in July 2023, and is an introduction to a category-theoretic perspective on metric spaces.  A metric space is a set of points such that between each pair of points there is a number -- the distance -- such that the triangle inequality is satisfied; a small category is a set of objects such that between each pair of objects there is a set -- the hom-set -- such that elements of the hom-sets can be composed. The analogy between the structures that can be made in to a common generalization of the two structures, so that both are examples of enriched categories.   This gives a bridge between category theory and metric space theory.  I will describe this and three examples from around mathematics where this perspective has been useful or interesting.  The examples are related to the tight span, the magnitude and the Legendre-Fenchel transform.
\end{abstract}

\setcounter{tocdepth}{2}
\tableofcontents

\section{Introduction}
This paper is an expanded version of the talk I gave at the CATMI meeting in Bergen, July 2023.  It is aimed at a wide audience across mathematicians and computer scientists, the main prerequisite is an interest in category theory, but no deep knowledge of category theory is required.  The purpose of the talk was to illustrate how category theory gives an interesting and fruitful perspective on metric spaces and related notions, not by considering the category of metric spaces, but rather by considering metric spaces (and certain generalizations) as being category-like entities themselves.  This is similar to the perspective of groups and monoids as one-object categories and posets as categories with at most one morphism between each pair of objects.

I will illustrate the perspective with three vignettes of mathematics that I've been directly involved in.  In each of the three I will describe some mathematical structure and then explain what the categorical perspective brings.  Any one of these might be viewed as a curiosity in itself, but taken together these demonstrate the broad perspective and show that category-like structures are likely to remain a rich vein with much still to be mined.

\subsection{Historical context}
Before moving on to an overview of the vignettes, I will make some historical comments.  Fr\'echet introduced the abstract notion of metric space under the name ``classe (E)''%
\footnote{``E-class'' would be a reasonable translation} in his 1906 thesis~\cite[p.~30]{Frechet:Thesis}.  This abstract treatment of sets with a notion of distance seems to owe much to the abstract treatment of group theory in the late nineteenth century.  (See Taylor's work~\cite{Taylor_1982} on Fr\'echet's thesis for more background.)  This definition was picked up by Hausdorff in his influential 1914 book~\cite[p.~211]{Hausdorff:Grundzuge} on sets and general topology, and was given the familiar presentation and name ``metrischer Raum'' (``metric space'' in English).  The key axiom for a metric space is the triangle inequality:
\begin{equation}
    \dd(x, x') + \dd(x', x'') \ge \dd(x, x'').
    \label{eq:triangle-ineq}
\end{equation}

Eilenberg and Mac Lane introduced the notion of category in 1945~\cite{EilenbergMacLane:General-theory} and a key part of the structure of a category is the composition of morphisms which can be axiomatized as a function on the product of hom-sets:
\begin{equation}
    \Hom(c, c') \times \Hom(c', c'') \to \Hom(c, c'').
    \label{eq:composition}
\end{equation}

Whilst attending a seminar in 1967, Lawvere observed the formal similarity between~\eqref{eq:triangle-ineq} and~\eqref{eq:composition} and this laid the basis for his important 1973 paper~\cite{Lawvere:metric-spaces} on metric spaces as enriched categories; enriched categories having been introduced around 1965 by Eilenberg and Kelly~\cite{EilenbergKelly:Closed-categories} and also by Linton (see Street's historical comments~\cite[p.~238]{Street:conspectus}).  It is this work of Lawvere which is the foundation for the work represented here.  The necessary background is presented in Section~2 below.  You'll see how various bits of metric space theory fit naturally into a category theoretic context and how more general notions of `metric space' can sometimes be appropriate; for instance the Hausdorff metric on the subsets of a metric space is more naturally seen to be one of Lawvere's more general metric spaces (aka \(\Rplusbar\)-categories) where the distance is not symmetric.

\subsection{The three vignettes}
Here I'll briefly summarize the three vignettes considered in the paper.  

The first vignette involves the \define{tight span} of a metric space, this is a construction which goes by many names as it has been discovered many times and has many applications, in particular in the theory of network flow.  The tight span is a minimal nice metric space in which the original metric space embeds.  There are recent generalizations of this to asymmetric metric spaces (in particular for directed flow) and from the perspective of category theory this generalization can be seen as an instance of an Isbell, or profunctor, nucleus.  The classical tight span can be characterized in these terms as well.

The second vignette involves an invariant of metric spaces called the \define{magnitude}.  This was first discovered by two mathematical biologists, Solow and Polasky, under the name `effective number of species' whilst considering notions of biological diversity using a `utilitarian approach'.  Much later, Leinster was considering notions of `Euler characteristic' for finite categories and realised that it could be generalized to certain enriched settings, in particular to Lawvere's generalized metric spaces, this recovered the effective number of species and became known as magnitude.  Following this, magnitude was found to have a rich mathematical life, aside from the biological and enriched categorical connections it was found to have connections to homology theory, potential theory and classical notions of volume and dimension.

The third vignette involves the \define{Legendre-Fenchel transform} which appears in convex analysis, Hamiltonian physics and large deviation theory among other places.  For this we have a different generalization of the notion of classical metric space where we allow negative distances; this notion is obtained by enriching over a different base category, namely the whole real line (with plus and minus infinity) instead of just the non-negative half (with plus infinity).

\subsection{Why bother with category theory?}
A standard question to ask is ``Why bother with the abstraction of enriched category theory here''.  Here are a few responses to that question.
\begin{itemize}
    \item Sometimes it allows you to see how constructions are `categorically inevitable' -- to use Leinster's phrase~\cite{Leinster:Codensity}.
    \item It allows you to see how different seeming constructions are essentially the same.
    \item It can make different areas easier to understand if you know related constructions in a different area.
    \item It can allow for formal analogies to become mathematical theorems or mathematical functors.
    \item It can allow for generalizations.
    \item It can make clear what is natural in general (what is abstract nonsense) and what is specific to that context.
\end{itemize}

\subsection{Acknowledgements}
I would like to thank the organizers of the CATMI conferences for the opportunity to talk about this material.  I would also like to thank Bruce Bartlett and Jared Ongaro for inviting me to talk on a similar theme at the African Mathematical Seminar in 2022.  Thanks also to Callum Reader for some comments on a draft of this paper and to the anonymous referee for helpful feedback.

\section{Background on enriched category theory}
In this section we'll start with the general notions of enriched category, then show how metric-like structures fit into this framework and give some natural examples of generalized metric spaces.  We then go on to look at `scalar-valued functors' and the Yoneda embedding before finishing with the idea of the nucleus of a profunctor.

\subsection{Enriched categories}
\label{subsec:enriched-cats}

We will look at the idea of enriched categories here but not all the details.  For more details, I would recommend the seminal paper of Lawvere~\cite{Lawvere:metric-spaces}; for slightly more conventional treatments you might consider the chapter of Borceaux's book~\cite{Borceux:Handbook-two} or the section of Johnson and Yau's book~\cite{JohnsonYau:Two-dimensional-categories}.  More advanced treatments can be found in the books of Kelly~\cite{Kelly:Enriched-category-theory} and Riehl~\cite{Riehl:categorical-homotopy-theory}.  Full definitions can also be found at Wikipedia and the nLab~\cite{nlab:enriched_category}.

We will first recap the specific definition of category that we are going to generalize.  Actually the definition we will give is that of a ``locally small'' category, which means that the morphisms between a fixed pair of objects forms a set rather than a proper class.

\begin{defn}
    A \define{category} \(\C\) consists of a collection, \(\ob \C\), of \define{objects} together with the following data which are required to satisfy the unit and associativity axioms:
    \begin{enumerate}
        \item for all \(c, c'\in \ob \C\) there is a specified set \(C(c, c')\) called the \define{hom-set} from \(c\) to \(c'\);
        \item for all \(c, c', c''\in \ob \C\) there is a specified function called \define{composition},
        \[{\circ}\colon \C(c, c') \times \C(c', c'') \to \C(c, c'');\]
        \item for all \(c \in \ob \C\) there is a specified element \[\id_c\in \C(c, c).\]
    \end{enumerate}
\end{defn}

A starting point for generalizing this definition to that of enriched category is the observation that there are many categories where the hom-sets can be considered as `things' other than sets.

For example, consider the category \(\Vect\) of vector spaces and linear maps; for concreteness, let's work over the complex numbers \(\CC\).  For all vector spaces \(V, V'\in \ob \Vect\) the hom-set \(\Vect(V, V')\) has a canonical vector space structure -- inherited from that of \(V'\) -- and moreover each composition map of the form 
\[
    \Vect(V, V') \times \Vect(V', V'') \to \Vect(V, V'')
\]
is \emph{bilinear} and so corresponds to a linear map from the tensor product
\[
    \Vect(V, V') \otimes \Vect(V', V'') \to \Vect(V, V'').
\]
It is not important that you know what the tensor product of vector spaces is, but what is important here is the structure that it equips the category of vector spaces with, namely the structure of a monoidal category.  Then the right notion of `thing' that replaces `set' in hom-set is that of object in a monoidal category.

A \define{monoidal category} is a category \(\V\) equipped with a monoidal product which is a functor \({\otimes}\colon \V \times \V \to \V\), and a unit object \(\one \in \ob \V\), these are required to be associative and unital in an appropriate sense, but that does not need to bother us here.  (For more details see the references mentioned at the beginning of this subsection or the chapter in Mac Lane's book~\cite{MacLane:Categories-for-the-working}.)  We will write \((\V, {\otimes}, \one)\) for such a monoidal category.

We have two examples of monoidal categories above: \((\Vect, {\otimes}, \CC)\), the category of vector spaces equipped with the tensor product, having the ground field \(\CC\) for unit object; and \((\Set, {\times}, \terminalset)\) the category of sets with equipped with the cartesian product, having a one object set for the unit object.  

We can rewrite the definition of locally small category to emphasize where the monoidal structure of the category of sets has been used, so that we can easily generalize the definition.

\begin{defn}
    A \define{category} \(\C\) consists of a collection, \(\ob \C\), of \define{objects} together with the following data which are required to satisfy the unit and associativity axioms:
    \begin{enumerate}
        \item for all \(c, c'\in \ob \C\) there is a specified object in the category of sets, \(C(c, c')\in \ob\Set\), called the \define{hom-set};
        \item for all \(c, c', c''\in \ob \C\) there is a specified morphism in \(\Set\) called \define{composition},
        \[{\circ}\colon \C(c, c') \times \C(c', c'') \to \C(c, c'');\]
        \item for all \(c \in \ob \C\) there is a specified morphism in \(\Set\), \[\id_c \colon \terminalset \to \C(c, c).\]
    \end{enumerate}
\end{defn}

It is then a small jump to the following generalization.
\begin{defn}
    For \((\V, {\otimes}, \one)\) a monoidal category, a \define{category enriched over~\((\V, {\otimes}, \one)\)}, or, simply, a \define{\(\V\)-category}, \(\C\), consists of a collection, \(\ob \C\), of \define{objects} together with the following data which are required to satisfy the unit and associativity axioms (see the references above for details):
    \begin{enumerate}
        \item for all \(c, c'\in \ob \C\) there is a specified object, \(C(c, c')\in \ob\V\) called the \define{hom-object};
        \item for all \(c, c', c''\in \ob \C\) there is a specified morphism in \(\V\) called \define{composition},
        \[{\circ}\colon \C(c, c') \otimes \C(c', c'') \to \C(c, c'');\]
        \item for all \(c \in \ob \C\) there is a specified morphism in \(\V\), \[\id_c \colon \one \to \C(c, c).\]
    \end{enumerate}
\end{defn}

It should be clear that a category enriched over \((\Set, \times, \terminalset)\) is precisely a locally small category.  Also, there is a \(\Vect\)-enriched category -- let's call it \(\selfenr{\Vect}\) --
associated to the category \(\Vect\); the objects are vector spaces and the hom-object \(\selfenr{\Vect}(V, V')\) is the \emph{vector space} of linear maps from \(V\) to \(V'\).  Of course, for \(\Vect\) the comparable hom-set \(\Vect(V, V')\) is the \emph{set} of linear maps from \(V\) to \(V'\).  The distinction between \(\Vect\) and \(\selfenr{\Vect}\) is subtle but important.  The subtlety in this case is mainly because \(\Vect\) is a category whose objects are sets with extra structure.  In order to think of metric-like spaces in terms of enriched category theory, we will need to look at  enriching over monoidal categories whose objects are not sets with extra structure, meaning that the enriched categories themselves do not have sets of morphisms between their objects, only hom-objects.

Before doing that, however, I should mention the generalization of functor in this setting.

\begin{defn}
    For \((\V, {\otimes}, \one)\) a monoidal category, and \(\V\)-categories \(\C\) and~\(\D\), a \(\V\)-functor consists of a function \(F\colon \ob{\C} \to \ob{\D}\), together with, for each \(c, c' \in \ob{\C}\), a morphism in \(\V\) between hom-objects \(F\colon \C(c, c') \to \D(F(c), F(c'))\) and these are required to satisfy conditions saying they behave nicely with regard to the identities and to composition (see references above for details).
\end{defn}

\subsection{Metric spaces as enriched categories}

Following the insights of Lawvere~\cite{Lawvere:metric-spaces},  we are interested in enriching over the following monoidal category in order to look at metric spaces.

\begin{defn}
    The category \(\Rplusbar\) has as objects non-negative real numbers together with infinity, so \(\ob \Rplusbar = \R_{\ge 0} \cup \{\infty\}\).  There is a morphism from \(a\) to \(b\) precisely when \(a\ge b\) and no morphism otherwise.  The monoidal category \((\Rplusbar, {+}, 0)\) is obtained by equipping \(\Rplusbar\) with addition as the monoidal product, where anything added to infinity gives infinity.
\end{defn}

The objects of \(\Rplusbar\) are real numbers and I know of no helpful way to think of them as having elements, i.e., as sets.  This means that \(\Rplusbar\)-categories have a slightly different feel to ordinary categories as there are no morphisms between a given pair of objects, there is only a hom-object between the objects, which, in this case, is a number.   

Unpacking the definiton of an \(\Rplusbar\)-category we get the following.
\begin{prop}
    An \define{\(\Rplusbar\)-category}, \(X\), consists of a collection, \(\ob X\), of {objects} together with the following data:
    \begin{enumerate}
        \item for all \(x, x'\in \ob X\) there is a specified number, \(X(c, c')\in \R_{\ge 0} \cup \{\infty\}\);
        \item for all \(x, x', x''\in \ob X\) there is the inequality,
        \[X(x, x') + X(x', x'') \ge X(x, x'');\]
        \item for all \(x \in \ob X\) there is an equality \(0 = X(x, x)\).
    \end{enumerate}
    Associativity and unitality are automatically satisfied due to \(\Rplusbar\) being \emph{thin}, that is, having at most one morphism between each pair of objects.
\end{prop}

It is sensible to interpret the number \(X(x, x')\) as a distance from \(x\) to \(x'\) and the similarity to the classical notion of metric space \'a la Fr\'echet should be apparent.  There are some notable differences however.

\begin{enumerate}
    \item The distances are not necessarily symmetric, so you could have \(X(x, x') \ne X(x', x)\).  This is not unreasonable, certainly in real life getting from one place to another might require a different route on the return journey.
    \item Distances can be infinite.  An infinite distance from \(x\) to \(x'\) represents the impossibility of getting from \(x\) to \(x'\).
    \item The distance from one point to a different point can be zero.  We will see some examples below showing that this is reasonable.
\end{enumerate}

From this generalized metric space perspective, an \(\Rplusbar\)-functor \(f\colon X \to Y\) is viewed as a \define{short map} or distance non-increasing function, so it is a function \(f\colon \ob X \to \ob Y\) such that \(X(x, x') \ge Y\big(f(x), f(x')\big)\) for all \(x, x' \in \ob X\).

\subsection{Examples of \texorpdfstring{\(\Rplusbar\)}{R+}-categories}
Let's now look at some examples; in particular we will see natural examples where  Frech\'et's axioms are weakened.
\subsubsection{Self enrichment} 
Just as we did with \(\Vect\) above, we can create an \(\Rplusbar\)-category \(\selfenr{\Rplusbar}\) out of \(\Rplusbar\) itself, or at least its objects.  (This corresponds to \(\Rplusbar\) being a closed monoidal category.)  The objects of \(\selfenr{\Rplusbar}\) are the objects of \(\Rplusbar\), and the hom-objects, or generalized distances, are defined using the \define{truncated difference}, so for \(a, b\in [0, \infty]\) we have
\[
    \selfenr{\Rplusbar}(a, b) = b \minusdot a 
    \coloneq
    \max(b - a, 0),
\]
where \(\infty \minusdot \infty \coloneq 0\) (this can be contrasted with what you will see in the third vignette).
If we think of the extended non-negative real numbers as standing vertically and going upwards, then for \(b \ge a\) we can think of it being `free' to descend from \(b\) to \(a\), but having a `cost' of \(b-a\) associated with ascending from \(a\) to \(b\).
%
%
\[
    \raisebox{-.45\height}{\input{\figdir/Rplusbar_metric.pstex_t}}
\]

\subsubsection{Spaces of short maps}
\label{subsubsec:spaces-of-short-maps}
As the category \(\Rplusbar\) is sufficiently nice -- it has all categorical limits, which are given by suprema -- if \(X\) and \(Y\) are \(\Rplusbar\)-categories then there is the \define{functor \(\Rplusbar\)-category} \(\funcat{X}{Y}\) consisting of the short maps, or distance-decreasing maps, from \(X\) to \(Y\) and the generalized metric given as follows:
\[
    \funcat{X}{Y}(f, g) \coloneq \sup_{x\in X} Y(f(x), g(x)).
    \qquad
    \raisebox{-.45\height}{\input{\figdir/function_metric.pstex_t}}
\]
This is a measurement of the furthest apart that \(f\) and \(g\) get.  

In particular, if we take \(Y= \selfenr{\Rplusbar}\), which is \(\Rplusbar\) with the truncated difference metric as above, then we can get the following generalized metric on `scalar-valued' short maps: 
\begin{equation}
    \funcat{X}{\selfenr{\Rplusbar}}(f, g) \coloneq \sup_{x\in X} \left(g(x) \minusdot f(x)\right).
    \label{eq:scalar-valued-fns-metric}
\end{equation}
This can be viewed as a refinement of the standard `sup metric' on the space of bounded maps from a set to the non-negative real numbers.  If \(P\) is a set then it can be considered as an \(\Rplusbar\)-category with the \define{discrete} metric: \(P(p, p') = \infty\) if \(p \ne p'\).  This means that all functions \(P\to [0, \infty]\) can be thought of as short maps \(P\to \selfenr{\Rplusbar}\).  If \(f, g \colon P \to [0,\infty]\) are bounded functions then the usual sup metric is given by symmetrizing the above asymmetric metric:
\[
    \sup_{p\in P}\left| g(p) - f(p) \right| = \max\left(\funcat{P}{\selfenr{\Rplusbar}}(f, g), \funcat{P}{\selfenr{\Rplusbar}}(g, f) \right). 
\]



\subsubsection{Asymmetric Hausdorff metric on subsets of a metric space} 

Now we can see an example where even if we are just interested in classical metric spaces we naturally end up with an asymmetric metric.  In fact it is an asymmetric refinement of the well-known -- and symmetric -- Hausdorff  metric.  

Given a classical metric space \(M\) we can take the set of closed subsets, \(S_M \coloneq \{A \subseteq M \mid A\ \text{closed}\}\).  We can then define a generalized metric on this by
\[
    S_M(A, B) \coloneq \adjustlimits\sup_{a\in A}\inf_{b\in B} M(a, b).
    \qquad
    \raisebox{-.45\height}{\input{\figdir/compact_subspaces.pstex_t}}
\]
One way of thinking about \(S_M(A, B)\) -- at least for \(A\) and \(B\) non-empty -- is that it is measuring the furthest you would have to go if you were dropped at a random point in \(A\) and wanted to take the shortest route to \(B\).  Thinking about the empty set, we have have \(S_M(\emptyset, B) = 0\) for any subset \(B\) and \(S_M(A, \emptyset) = \infty\) for any non-empty subset \(A\).

Because we are considering \emph{closed} sets we have the following simple characterization of zero-distance, which is telling us that the generalized metric is encoding the usual order on subsets:
\[
    S_M(A, B) = 0 
    \quad \Longleftrightarrow \quad
    A \subseteq B.
\]
This leads us to see that this generalized metric is naturally asymmetric even though we started with a symmetric metric on \(M\).  The usual Hausdorff metric \(\hausdorff\) -- defined for \emph{non-empty}, \emph{compact} subsets -- is obtained by symmetrizing this generalized metric:
\[
    \hausdorff(A, B) = \max\big(S_M(A, B), S_M(B, A)\big).
\]
Clearly, however, the Hausdorff metric is losing information, such as the partial order, that was contained in the generalized metric.

\subsection{Scalar-valued functors and the Yoneda embedding}

One philosophy or approach to working with \(\V\)-enriched categories is to think of it as like working over a given field \(\F\) in algebra.  Certainly in the common, nice case of the existence of a \(\V\)-enriched version \(\selfenr{\V}\) of \(\V\) we can think of \(\selfenr{\V}\) or \(\V\) as being like a ground field or a collection of scalars.  In many areas of mathematics, it is fruitful to probe or study an entity such as a set, topological space, manifold, module, etc., by considering the totality of (possibly nicely-behaved) scalar-valued functions on that entity.

For example, for a set \(S\) one considers \(\Fun(S, \F)\) the vector space of \(\F\)-valued functions on \(S\), which might be written \(\F^S\) or \(\F[S]\).  The set \(S\) embeds in the space of functions \(\Fun(S, \F)\) via the `delta-function':
\begin{equation}
    S\hookrightarrow \Fun(S, R); 
    \quad 
    s \mapsto \Biggl( \delta_s\colon s'\mapsto 
    \begin{cases}
        1 & \text{if }s=s' \\
        0 & \text{otherwise}
    \end{cases}\,
    \Biggr).
    \label{eq:delta-function}
\end{equation}

One carries out a similar process in the context of \(\V\)-enriched categories.  If \(\V\) is suitably nice -- closed, symmetric monoidal -- such as the examples \(\Set\), \(\Vect\) and \(\Rplusbar\) used above, then there is a \(\V\)-enriched category \(\selfenr{\V}\) which is the \(\V\)-enriched version of \(\V\) itself.  Provided that \(\V\) has sufficiently many colimits (for example, is cocomplete), then for any small \(\V\)-category \(\C\) there are the following associated \(\V\)-categories of `scalar-valued' \(\V\)-functors:
\[
    \funcat{\C}{\selfenr{\V}}
    \quad\text{and}\quad
    \funcat{\C^\op}{\selfenr{\V}}.
\]

These are fundamental objects, analogous to the scalar valued functions on a set.  For historical reasons, these are sometimes referred to as the enriched categories of, respectively, enriched copresheaves and enriched presheaves on \(\C\), however, in general, this terminology is more scary than helpful and we will refer to them as the enriched categories of, respectively, covariant and contravariant scalar-valued, or \(\selfenr{\V}\)-valued enriched functors.

One of the fundamental results of enriched category theory is that the \define{Yoneda functor}, 
\[
    \C \to \funcat{\C^\op}{\selfenr{\V}}, \quad
    c \mapsto (c' \mapsto \C(c', c)),
\]  
is a fully faithful functor, meaning there it gives isomorphisms on hom-objects.
The Yoneda functor is the analogue of the delta-function~\eqref{eq:delta-function}.   Let's look at two examples.
\begin{enumerate}
    \item If we consider \(\Vect\)-categories and have \(A\) an algebra over \(\CC\) -- that is a vector space with an associative, unital product -- then there is the associated one-object \(\Vect\)-category \(\bar{A}\).  The \(\Vect\)-category of \(\selfenr{\Vect}\)-valued \(\Vect\)-functors \(\funcat{\bar{A}^\op}{\selfenr{\Vect}}\) is the \(\Vect\)-category of left representations of \(A\).  The Yoneda functor picks out the left regular representation of \(A\), that is to say the action of \(A\) on itself via left multiplication.

    \item Consider an \(\Rplusbar\)-category, or generalized metric space, \(X\).  Here the functor \(\Rplusbar\)-category \(\funcat{X^\op}{\selfenr{\Rplusbar}}\) is described in Subsection~\ref{subsubsec:spaces-of-short-maps} above.  The objects, or points, are the distance non-increasing maps or short maps, and the metric is given in~\eqref{eq:scalar-valued-fns-metric}.  The Yoneda functor in this case \(X \to \funcat{X^\op}{\selfenr{\Rplusbar}}\), \(x\mapsto X({-}, x)\), is an isometry and gives a generalization of the classical Kuratowski embedding theorem.
\end{enumerate}

\subsection{The nucleus of a profunctor}
A profunctor is a categorical analogue of a matrix.  Associated to each profunctor is a certain pair of adjoint functors; I will motivate this construction by showing an analogous construction of a pair of adjoint \emph{functions} associated to a matrix.  The nucleus of the profunctor is obtained from this adjunction, and the nucleus will appear in the first and third vignette below.  

The nucleus is probably the least standard piece of category theory in this paper and I do not know of any textbook treatments.  I discuss it in my paper~\cite{Willerton:legendre-fenchel}, having learnt of it from Pavlovic~\cite{Pavlovic_2012}; it had previously appeared in the work of Gutiérrez García, Mardones-Pérez, de Prada Vicente and Zhang~\cite{GutierrezGarcia:fuzzy-galois}.

Suppose that \(R\) and \(S\) are finite sets; in practice they might be ordinal sets, \(R=\{1, \dots, m\}\) and \(S=\{1, \dots, n\}\).  An \(R\)-\(S\)-indexed matrix is an \(\F\)-valued function \(M\colon R\times S\to \F\).

Given such a matrix, we get a pair of linear maps between the spaces of scalar-valued functions:
\[
    \begin{tikzcd}[ampersand replacement=\&]
        {\F^R} \& {\F^S}.
        \arrow["{\overline{M}^{\mathrm{T}}}", curve={height=-6pt}, from=1-1, to=1-2]
        \arrow["{\overline{M}}", curve={height=-6pt}, from=1-2, to=1-1]
    \end{tikzcd}
\]
These are given by the familiar multiplication of a vector by a matrix:
\[
    \Bigl(\overline{M}w\Bigr)(r) \coloneq \sum_{s\in S}M(r,s)\cdot w(s)\,;
    \quad
    \Bigl(\overline{M}^{\mathrm{T}}v\Bigr)(s) \coloneq \sum_{r\in R}M(r,s)\cdot v(r).
\]
Furthermore, these linear maps are adjoint with respect to the canonical inner products on the function spaces \(\F^R\) and \(\F^S\):
\[
    \innerprod[\Big]{\overline{M}w}{v}{\F^S} = 
    \innerprod[\Big]{w}{\overline{M}^{\mathrm{T}}v}{\F^R}.
\]

A matrix is a scalar-valued function on a product of sets; in the \(\V\)-category setting the analogue of a matrix is a \define{profunctor} which is a scalar-valued \(\V\)-functor on a product of \(\V\)-categories, contravariant in one of the variables.  So if \(\C\) and \(\D\) are \(\V\)-categories then a profunctor, denoted \(P\colon \C \profto \D\), is a \(\V\)-functor of the form
\[
    P\colon \C^\op \otimes \D \to \selfenr{\V}.
\]
For example, for a \(\V\)-category \(\C\) the hom-objects assemble into a \(\V\)-functor,
\(
    \C({-},{-}) \colon \C^\op \otimes \C \to \selfenr{\V},
\)
called the \define{hom-profunctor}.

Given a profunctor as above, assuming that \(\V\) is sufficiently nice whilst \(\C\) and \(\D\) are sufficiently small, we get a pair of \(\V\)-functors between the \(\V\)-cateogries of scalar-valued \(\V\)-functors:
\begin{equation}
    \begin{tikzcd}[ampersand replacement=\&]
        {\funcat{\C^\op}{\selfenr{\V}}} \& {\funcat{\D}{\selfenr{\V}}^\op}
        \arrow["{P^\ast}",shift left, curve={height=-2pt}, from=1-1, to=1-2]
        \arrow["{P_\ast}", shift left, curve={height=-2pt}, from=1-2, to=1-1].
    \end{tikzcd}
    \label{eq:isbell-type-adjunction}
\end{equation}
These are given by `\(\V\)-ends' in a form which looks akin to matrix multiplication:
\[
    P_\ast(\phi)(c) = \int_{d\in \D} [\phi(d), P(c, d)];
    \qquad
    P^\ast(\psi)(d) = \int_{c\in \C} [\psi(c), P(c, d)].
\]
It is not easy to give a quick idea of what a \(\V\)-end is, so I won't go into that here -- see the references at the beginning of Section~\ref{subsec:enriched-cats}.  Hopefully it suffices to say that when \(\V\) is \(\Rplusbar\) these ends have the following form:
\[
    P_\ast(\phi)(c) = \sup_{d\in \D} \left\{ P(c, d) \minusdot \phi(d)\right\};
    \qquad
    P^\ast(\psi)(d) = \sup_{c\in \C} \left\{P(c, d) \minusdot \psi(c)\right\}.
\]

In general, \(P^\ast\) and \(P_\ast\) are adjoint \(\V\)-functors which and so we have natural isomorphisms of \(\V\)-objects:
\[
    \funcat{\D}{\selfenr{\V}}^\op\left(P^\ast(\psi), \phi\right)
    \cong
    \funcat{\C^\op}{\selfenr{\V}}\left(\psi, P_\ast(\phi)\right).
\]
The similarity to the matrix example above should be striking.  We will refer to the adjunction \eqref{eq:isbell-type-adjunction} as the \define{Isbell type adjunction} coming from the profunctor \(P\).

The \define{nucleus} of the profunctor \(P\) is then defined to be the \define{centre} of this adjunction.  This can be taken to be \(\Fix (P_\ast P^\ast)\) the fixed point category of \(P_\ast P^\ast\), that is the full subcategory of \(\funcat{\C^\op}{\selfenr{\V}}\) whose objects are functors \(\psi\) such that the `unit' \(\one \to \funcat{\C^\op}{\selfenr{\V}}(\psi, P_\ast P^\ast(\psi))\) from the adjunction is suitably invertible.  This can equivalently by taken to be \(\Fix(P^\ast P_\ast)\), and we get an adjoint equivalence by restricting the profunctor adjunction:
\[
\begin{tikzcd}[ampersand replacement=\&, column sep=tiny]
	{\Fix(P_\ast P^\ast)} \& \cong \& {\Fix(P^\ast P_\ast).}
	\arrow["P^\ast", shift left, curve={height=-2pt}, from=1-1, to=1-3]
	\arrow["{P_\ast}", shift left, curve={height=-2pt}, from=1-3, to=1-1]
\end{tikzcd}
\]
The profunctor nucleus will appear in both the first and third vignettes.

\section{First vignette: the tight span}
In the first subsection here we will look at the tight span from a classical perspective and in the second subsection we will look at it from the enriched category theoretic perspective.

\subsection{The classical tight span}

Suppose that \(M\) is a classical metric space, then there is a certain associated classical metric space \(T(M) \subset \Fun\left(M, [0, \infty)\right)\).  This has been discovered many times (for instance, by Isbell~\cite{Isbell:SixTheoremsInjective}, Dress~\cite{Dress:TreesTightExtensions} and Chrobak and Larmore~\cite{ChrobakLarmore}) and it goes by many names such as hyperconvex hull, injective envelope, Isbell hull and tight span.  We will call it the \define{tight span}.

The are various ways to define the tight span, but one definition is that it consists of all functions \(f \colon M \to [0, \infty)\) such that
\begin{enumerate}
    \item \(M(m, m')\ge f(m') - f(m)\) for all \(m, m' \in M\), 
    
    \item \(f(m) = \sup_{m'\in M}\{M(m, m')- f(m')\}\) for all \(m \in M\).
\end{enumerate}

The Kuratowski embedding gives an isometric embedding \(M \hookrightarrow T(M)\).  An example of a three-point space is given in Figure~\ref{fig:embedding-in-tight-span} where the tight span is a tripod with the lengths of the legs indicated.


\begin{figure}[ht]
    \centering
    \(M=~
    \begin{tikzpicture}[scale=0.7, baseline=-0.8cm]
        \node [circle,draw=blue,fill=blue,thick, inner sep=0pt,minimum size=1mm,label=above:$b$](x) at (0,0) {};
        \node [circle,draw=blue,fill=blue,thick, inner sep=0pt,minimum size=1mm,label=above:$c$](y) at (3,0) {};
        \node [circle,draw=blue,fill=blue,thick, inner sep=0pt,minimum size=1mm ,label=below:$a$](z) at (1,-1.7) {};
        \begin{scope}[<->,shorten >=2mm,shorten <=2mm]
            \path (x) edge node[above] {$r$} (y);
            \path (z) edge node[pos=0.5,left] {$t$} (x);
            \draw (y) edge node[pos=0.5,right] {$s$} (z);
        \end{scope}
    \end{tikzpicture}
    \quad\hookrightarrow\quad 
    T(M)=~
    \begin{tikzpicture}[scale=0.7, baseline=-0.8cm]
        \node [circle,draw=blue,fill=blue,thick, inner sep=0pt,minimum size=1mm,label=above:$b$](x) at (0,0) {};
        \node [circle,draw=blue,fill=blue,thick, inner sep=0pt,minimum size=1mm ,label=above:$c$](y) at (3,0) {};
        \node [circle,draw=blue,fill=blue,thick, inner sep=0pt,minimum size=1mm,label=below:$a$](z) at (1,-1.71) {};
        \node [circle,draw=blue,fill=blue,thick, inner sep=0pt,minimum size=0mm ]
        (o) at (1,-0.866) {};
        \draw [-,draw=blue,very thick] (o) 
            edge node[pos=0.4,left] {$\tfrac{r+t-s}{2}$\;} (x) 
            edge  node[pos=0.5,above] {$\tfrac{r+s-t}{2}$} (y)
            edge  node[midway,right] {$\tfrac{s+t-r}{2}$} (z);
    \end{tikzpicture}
    \)
    \caption{}
    \label{fig:embedding-in-tight-span}
\end{figure}

Generally, the tight span \(T(M)\) is the smallest `hyperconvex', or equivalently `injective' metric space in which \(M\) isometrically embeds.

There have been many applications of this construction; here are three.

\begin{description}
    \item [Recreating phylogenetic trees~\cite{Dress:TreesTightExtensions}.] If you have a set of species and some notion of distance between them then you can try to use this to recreate the phylogenetic tree which encodes the evolutionary history of the species.  The tight span can be used as a first approximation to the phylogenetic tree, indeed a metric space embeds in a tree if and only if the tight span is a tree.

    \item[Placing servers on a network~\cite{ChrobakLarmore}.] The \(k\)-server problem is about moving \(k\) around a metric space in an optimal fashion.  Chrobak and Larmore utilized the tight span to prove a conjecture on this when the network is a tree.
    
    \item[Multicommodity flow~\cite{Chepoi:TX-approach}.] One multicommodity flow problem is where you have a network -- a graph -- with some capacity restriction on the edges, you also have a set of pairs of nodes -- sources and sinks -- together with, for each pair, a demand level of commodities that must be transported from source to sink; the problem is to plan which edges to use to transport the commodities so that the demands are met but the capacities are not exceeded.  This is a linear programming problem and the dual problem involves cost metric on the edges; Chepoi applied tight span techniques to the dual problem.
\end{description}

\subsection{Enriched categories and the directed tight span}
(Material in this section is taken from my paper~\cite{Willerton:tight-spans}.)
Suppose that \(X\) is an \(\Rplusbar\)-category, so it could be a classical metric space.  We can consider the distance function as an \(\Rplusbar\)-profunctor, the \define{hom profunctor}:
\[
    X({-}, {-})\colon X^\op \times X \to \selfenr{\Rplusbar}.
\]
The Isbell-type adjunction \eqref{eq:isbell-type-adjunction} coming from the hom profunctor is the so-called \define{Isbell adjunction}
\[
\begin{tikzcd}[ampersand replacement=\&]
	{\funcat[\big]{X^\op}{\selfenr{\Rplusbar}}} \& {\funcat[\big]{X}{\selfenr{\Rplusbar}}^\op.}
	\arrow["", shift left, curve={height=-2pt}, from=1-1, to=1-2]
	\arrow["{}", shift left, curve={height=-2pt}, from=1-2, to=1-1]
\end{tikzcd}
\]
The nucleus of this adjunction is the \define{Isbell completion} \(I(X)\) which is an \(\Rplusbar\)-category.  \begin{figure}[thbp]
    \centering
    \begin{tikzpicture}[scale=0.3,line join=round]
        \draw[style =thin,opacity=1,arrows=->](0,0)--(2*5.033,2*1.355) node[below ] {$f(b)$};
        \draw[style =thin,arrows=->](0,0)--(1.3*7.27,1.3*-.938) node[below] {$f(a)$};
        
        \draw[style =thin,arrows=->](0,0)--(0,8.687) node[left] {$f(c)$};
        \filldraw[draw=black,fill=green!80,fill opacity=0.5](0,0)--(.808,-.104)--(.808,1.826)--(0,1.931)--cycle;
        \filldraw[draw=black,fill=green!80,fill opacity=0.5](.808,1.826)--(7.643,6.884)--(2.796,7.51)--(0,1.931)--cycle;
        \filldraw[draw=black,fill=green!80,fill opacity=0.5](.808,1.826)--(7.643,6.884)--(9.01,6.931)--(9.01,.174)--(.808,-.104)--cycle;
        
        \draw[draw=blue,style =ultra thick](2.796,7.51)--(3.542,3.849);
        
        \filldraw[draw=black,draw opacity=1,fill=green!80,fill opacity=0.5](.808,1.826)--(6.276,5.873)--(4.039,5.27)--cycle;
        
        \draw[draw=blue,style =ultra thick](4.039,5.27)--(3.542,3.849);
        \draw[draw=blue,style =ultra thick](9.01,.174)--(3.542,3.849);

        \draw [draw=blue,fill=blue] (4.039,5.27) circle (.15cm); 
        \node at (4.039,5.27) [anchor=south] {$a$};
        \draw [draw=blue,fill=blue] (2.796,7.51) circle (.15cm); 
        \node at (2.796,7.51) [anchor=south west] {$b$};
        \draw [draw=blue,fill=blue] (9.01,.174) circle (.15cm); 
        \node at (9.01,.174) [anchor=west] {$c$};
        \node [circle,draw=blue,fill=blue,thick, inner sep=0pt,minimum size=1mm,
        ](T) at (9.01,6.931){};
        \node [circle,draw=blue,fill=blue,thick, inner sep=0pt,minimum size=1mm,
        ](B) at (0,0) {};
    \end{tikzpicture}
    \caption{}
    \label{fig:directed-tight-span}
\end{figure}

The Yoneda embeddings are compatible with the Isbell adjunction and so give rise to an isometric embedding \(X\hookrightarrow I(X)\).  This embedding is pictured in Figure~\ref{fig:directed-tight-span} for the example of the three-point classical metric space from Figure~\ref{fig:embedding-in-tight-span}; the Isbell completion \(I(X)\) is shown as a subset of \(\R_{+}^{\{a, b, c\}}\), the asymmetric metric on this is given by 
\[
    I(X)(f, f') \coloneq \max_{x\in \{a, b, c\}}(f'(x)\minusdot f(x)).
\]

You can also see in Figure~\ref{fig:directed-tight-span} that the classical tight span has been drawn inside the Isbell completion and this is a general phenomenon for the Isbell completion of classical metric spaces as expressed in the following theorem.
\begin{thm}
    For \(M\) a classical metric space, the tight span \(T(M)\) sits isometrically inside the Isbell completion \(I(M)\) as the largest subset which contains \(M\) and for which the restriction of the generalized metric is classical, i.e., symmetric.
\end{thm}

From the category theory it turns out that the Isbell completion \(I(X)\) is both complete and cocomplete in an appropriate enriched sense and this translates to \(I(X)\) being a semi-module for the ``semi-tropical'' semi-ring \(([0, \infty], +, \max)\) in two different ways, giving a connection with tropical algebra.

At a similar time to the categorical investigations here, two other groups gave analogues of the tight span for generalized metric spaces, though not quite as generalized as \(\Rplusbar\)-categories as infinite distances were not considered.  In both cases the definition coincides with that of the Isbell completion.  
\begin{itemize}
    \item Hirai and Koichi defined~\cite{HiraiKoichi:tight-spans-directed-distances} a `directed tight span' motivated by mulitcommodity flows on directed networks.  
    \item Kemajou, Künzi and Olela Otafudu defined~\cite{Kemajou:isbell-hull-di-space} an `Isbell hull of a di-space' motivated by generalizing hyperconvexity and injectivity to directed metric spaces.
\end{itemize}

\section{Second vignette: the magnitude and the Euler characteristic}
In the first subsection we will see how mathematical biologists came up with essentially what is now called magnitude as an `effective number of species'.  In the second subsection we will look at how this was later obtained via an enriched category route, see some examples and comment on some connections with wider mathematics.

\subsection{Effective number of species}
Measuring biodiversity -- or more generally measuring diversity -- is about starting with a population, community or ecosystem of things (possibly organisms) and assigning some number representing the size of the diversity of the population.  Understanding what a good measure of diversity is is a problem for ecologists and geneticists in particular.  Magurran's book~\cite{Magurran:Measuring-bio-diversity} discusses this in broad terms.  

There are various desiderata that you might wish of such a measure and these might vary depending on the context; different properties might be of interest in different situations.  For instance, a standard, traditional way to mathematically model a population is as a finite probability space.  Each point corresponds to a species in the population and the probability assigned to a point represents the relative abundance of the corresponding species in the population; how the abundance is quantified, e.g., by number of individuals or total mass, is a  question which we won't consider here.  Many diversity measures are ways of getting numbers from finite probability spaces and quite a few are related to entropy.  One family of such diversity measures is the family of so-called Hill measures which we will see in the subsection on diversity measures below.

Solow and Polasky~\cite{Solow-Polasky:measuring-bio-diveristy} took a ``utilitarian approach''; they wanted a measure that would indicate the potential utility of a population, for example, in providing a future cure for a disease.  Here the relative abundances of the species are not as important as the differences between the species, so they mathematically modelled a population as a metric space \((S, d)\): points again correspond to species in the population, and the distance is some measure of difference between species.  Again, how the difference is actually quantified, e.g., genomically or physiologically, is a question we won't consider here.  Solow and Polasky also require a decreasing function \(f \colon [0, \infty] \to (0,1)\) which can be thought of as converting `species difference' to `species similarity'.  As an example they picked \(f(x)\coloneq e^{-\theta x}\) for some fixed \(\theta > 0\).  They defined the `similarity matrix' \(F=\{F_{i,j}\}_{i,j\in S}\) by \(F_{i,j}=f(d(i, j))\) and were lead to define the diversity measure \(V(S)\), in the case that the similarity matrix \(F\) is invertible, by
\begin{equation}
    V(S)= \sum_{i,j\in S}F^{-1}_{i,j}.
    \label{eq:effective-number}
\end{equation}
They studied some examples and basic properties of this measure which they showed could be thought of as an ``effective number of species''.  Now we'll look at a rather different way to approach this.

\subsection{Magnitude and Euler characteristic}
There are many notions of `size' of things in mathematics.  One recurrent motif is that of Euler characteristic.  
Euler characteristic is perhaps most well known as an invariant of surfaces, but more generally it is defined for topological spaces which have suitable finiteness properties.  However, it is also an important notion in algebra and combinatorics as well.

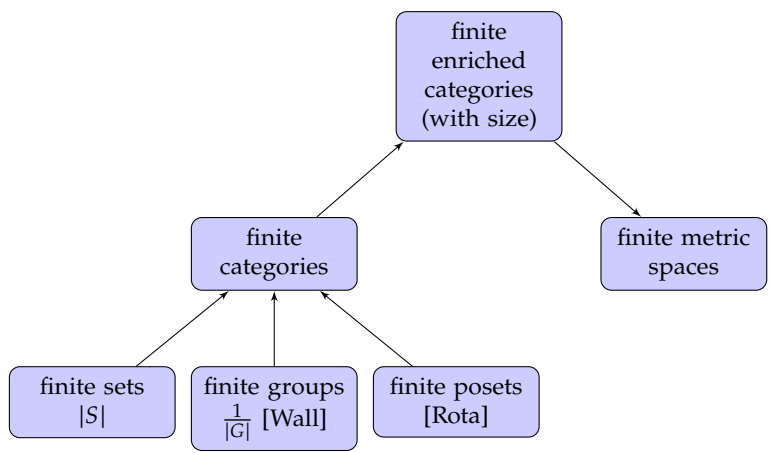
\begin{figure}[ht]
    \centering{
        \begin{tikzpicture}[font=\footnotesize]
            \node [block] (encategory) {finite enriched categories (with size)};
            \node [block, below left=1 and 0.5 of encategory] (category) {finite categories};
            \node [block, below right=1 and 0.5 of encategory] (metric) {finite metric spaces};
            \node [block, below right= 1 and 0.2 of category] (poset) {finite posets\\{[Rota]}};
            \node [block, below left= 1 and 0.2 of category] (set) {finite sets\\$\card{S}$};
            \node [block, below= 1 of category] (group) {finite groups\\$\frac{1}{\card{G}}$ [Wall]};
            \path [line] (poset) -- (category);
            \path [line] (set) -- (category);
            \path [line] (group) -- (category);
            \path [line] (category) -- (encategory);
            \path [line] (encategory) -- (metric);
        \end{tikzpicture}
    }
    \caption{A schematic of how Leinster generalized several notions of Euler characteristic and then specialized to finite metric spaces}
    \label{fig:schematic-of-generalization}
\end{figure}

If you have a finite set then its Euler characteristic is precisely the cardinality of the set.  For a finite group, Wall~\cite{Wall:Euler-char-group} showed that the right notion of Euler characteristic is the reciprocal of the order of the group.  For a finite poset, Rota~\cite{Rota:poset} gave a definition of the Euler characteristic.  Leinster~\cite{Leinster:Euler-char-category} observed that all of these three mathematical entities can be viewed as finite categories and he gave a definition of Euler characteristic of a finite category which generalized these three examples.   Let's have a look at what this is.

For a finite category \(\C\), define \(Z\), the \define{similarity matrix} of \(\C\), to be the \(\ob \C\)-\(\ob\C\)-indexed matrix such that \(Z(c, c') \coloneq \card{\C(c, c')}\) is the cardinality of the hom-set.  In the case that the similarity matrix \(Z\) is invertible, we can define the Euler characteristic or \define{magnitude} of \(\C\) as the sum of the entries of the inverse of \(Z\):
\[
    \magnitude{\C} \coloneq \sum_{c,c'\in \C}Z^{-1}(c,c') \in \Q.
\]
If the similarity matrix is not invertible then sometimes it is still possible to define the magnitude, but we won't need the details here.

Leinster then observed~\cite{Leinster:magnitude} that this definition can be generalized to a notion of Euler characteristic, or magnitude, for finite \(\V\)-categories provided that the enriching monoidal category \((\V, \otimes, \one)\) is equipped with a \define{size map} taking values in a commutative, semiring \(k\): a size map is a function \(\size{{\cdot}}\colon \ob\V\to k\) such that  
\begin{enumerate}[label=(\roman*)]
    \item \(\size{v} = \size{v'}\) if \(v\cong v'\); 
    \item \(\size{v \otimes v'} = \size{v}\cdot \size{v'}\); and 
    \item \(\size{\one} = 1\).
\end{enumerate}
In the case of the category of finite sets, cardinality gives a size map, so there is no contradiction with the notation used above.

If \((\V, \otimes, \one)\) is equipped with a {size map} then for a \(\V\)-category \(\C\) you can define the \define{similarity matrix} \(Z\) by \(Z(c, c') \coloneq \size{\C(c, c')}\) analogous to the \(\Set\)-category case done above and the magnitude of \(\C\) can then be defined from the similarity matrix \(Z\) as above.   (See Figure~\ref{fig:schematic-of-generalization} for a schematic representation of the generalizations.)

For enrichment using finite dimensional vector spaces, \(\V=\FinVect\), we can take the size to be the dimension, \(\size{V} \coloneq \dim(V)\).  For the case of primary interest, \(\V= \Rplusbar\), multiplicativity of the size map means that addition gets sent to multiplication and a sensible choice of size map turns out to be \(a \mapsto e^{-a} \in \R\).

This means that for an \(\Rplusbar\)-category \(X\) the similarity matrix is given by \(Z(x, x') \coloneq e^{-X(x, x')}\).  Provided that this similarity matrix is invertible, the magnitude of \(X\) can be defined by 
\[
    \magnitude{X} \coloneq \sum_{x,x'\in X}Z^{-1}(x,x') \in \R.
\]

For a given metric space or \(\Rplusbar\)-category \(X\) it is fruitful to consider the magnitude not just of \(X\) itself, but of all the rescalings of \(X\).  For \(t\in [0, \infty)\) define \(tX\) to have the same objects as \(X\) but with the distances scaled by a factor of \(t\), so \((tX)(x, x') \coloneq t\cdot X(x, x')\).  We then define the \define{magnitude function} of \(X\) to be the function \(t \mapsto \magnitude{tX}\), though it might not be defined for all values of \(t\).  This is equivalent to replacing the size map \(a\mapsto e^{-a}\) with \(a \mapsto e^{-ta}\); the resulting magnitude function is precisely the same as Solow and Polasky's effective number of species~\eqref{eq:effective-number}, with the parameter \(\theta\) being called \(t\) here.  The connection was discovered after comments by a member of the audience in a seminar given by Leinster.

\begin{figure}[tbh]
    \centering
    \begin{subfigure}{\textwidth}
        \centering
        \[
            tX:=\quad
            \begin{tikzpicture}[node distance = 4cm, auto,baseline=-0.1em]
                \node [point] (one) at (0,0) {};
                \node [point] (two) at (6,0.2) {};
                \node [point] (three) at (6,-0.2) {};
                \path [line,<->] (one) -- node[above,pos=0.5] {t}  (two);
                \path [line,<->] (one) -- node[below,pos=0.5] {t}  (three);
                \path [line,<->] (three) -- node[right,pos=0.5] {0.001t}  (two);
            \end{tikzpicture}
        \]
        \begin{tikzpicture}
            \begin{axis}[
                    width = 0.8\textwidth,
                    axis x line=center, 
                    axis y line=left,
                    xmin=-5.2, ymin=0, 
                    ymax=3.2, xmax=2.5,
                    xtick={-5,...,2}, 
                    xticklabels={0.01, , 1, ,100, ,10000 ,},
                    ytick={0,1,2,3}, 
                    yticklabels={0,1,2,3},
                    x axis line style={style = -}, 
                    y axis line style={style = -},
                    xlabel=$t$, 
                    ylabel=$\magnitude{tX}$,
                    yscale=0.5,
                    ]
                \addplot[mark=none] file {magnitude_of_3_points.dat};
                \addplot[mark=none,red,dashed] expression {3};
            \end{axis}
        \end{tikzpicture}
        \caption{}
        \label{subfig:three-point-metric-space}
    \end{subfigure}
    \begin{subfigure}{\textwidth}
        \centering
        \(tY \coloneq~
        \begin{tikzpicture}[
            node distance=4cm, 
            auto, 
            baseline=5em,
            ]
            \node [point] (A) at (0,0) {};
            \node [point] (B) at (0,2) {};
            \node [point] (C) at (0,4) {};
            \node [point] (D) at (2,1) {};
            \node [point] (E) at (2,3) {};
            \path [line,<->] (A) -- node[above,pos=0.4] {t}  (D);
            \path [line,<->] (A) -- node[above,pos=0.3] {t}  (E);
            \path [line,<->] (B) -- node[above,pos=0.4] {t}  (D);
            \path [line,<->] (B) -- node[above,pos=0.3] {t}  (E);
            \path [line,<->] (C) -- node[above,pos=0.4] {t}  (D);
            \path [line,<->] (C) -- node[above,pos=0.3] {t}  (E);
        \end{tikzpicture}\)
        \qquad
        \begin{tikzpicture}[baseline=5em]
            \begin{axis}[
                width = 0.5\textwidth,
                axis x line=center, 
                axis y line = left,
                xmin=-2.2, ymin=-0.2, 
                ymax=5.5, xmax=2,
                xtick={-2,...,2}, 
                xticklabels={0.01, 0.1, $1$, 10, 100},
                ytick={0,1,2,3,4,5}, 
                yticklabels={0,1,2,3,4,5},
                x axis line style={style = -},
                y axis line style={style = -},
                xlabel=$t$, 
                ylabel=$\magnitude{tY}$
                ]
                \addplot[mark=none] file {BlockOfDotsMagBadFivePointsI.dat};
                \addplot[mark=none] file {BlockOfDotsMagBadFivePointsII.dat};
                \addplot[mark=none,red,dashed] expression {5};
            \end{axis}
        \end{tikzpicture}
        \caption{}
        \label{subfig:bad-metric-space}
    \end{subfigure}
    \caption{}
    \label{fig:magnitude-function-plots}
\end{figure}

Figure~\ref{fig:magnitude-function-plots} gives the graphs of the magnitude functions of two classical metric spaces.  
In Figure~\ref{subfig:three-point-metric-space} we see a three-point metric space \(X\) with two points much closer to each other than they are to a third.  As far as the magnitude is concerned there are roughly three scaling regimes.  When \(0 < t \ll 1\) then \(tX\) looks very much like a single point and has magnitude close to \(1\); when \(10 < t < 100\) then \(tX\) looks like two reasonably well separated points and the magnitude is close to \(2\).  When \(t \gg 0\) then \(tX\) looks like three well separated points and the magnitude is close to \(3\).  This example illustrates why it can be useful to think of the magnitude as being a measure of ``effective number of points''.  It also illustrates the following theorem.

\begin{thm}
    For \(X\) an \(\Rplusbar\)-category, if \(t \gg 0\) then the magnitude function \(t\mapsto\magnitude{tX}\) is increasing, and as \(t \to \infty\) the magnitude \(\magnitude{tX}\) tends to the number of points in \(X\).
\end{thm}

Figure~\ref{subfig:bad-metric-space} illustrates a less well-behaved example.  Here the magnitude function is not defined when \(t=\ln 2\).

Many metric spaces are better behaved than this.
If \(\{e^{-\dd(i,j)}\}_{i,j}\) is positive definite then \(|X|\) is defined.
For example, if \(X\) is a subset of Euclidean space then \(|X|\) is defined.  Notably, Meckes~\cite{Meckes:positive-definite-spaces} proved that the definition of magnitude could be extended to infinite metric spaces provided that they are compact are positive definite in the sense that the similarity matrix is positive definite for a finite sub-metric space.

\begin{figure}[tbhp]
    \begin{center}
        \begin{tikzpicture}[
            small mindmap, concept color=blue!30,
            level 1 concept/.append style={concept color=blue!15},
            ]
            \node [concept] {magnitude}
            child[grow=0] {node[concept] {Minkowski dimension}}
            child[grow=60 ] {node[concept] {potential theory and special functions}}
            child[grow=120 ] {node[concept] {magnitude homology}}
            child[grow=180 ] {node[concept] {enriched category theory}}
            child[grow=240 ] {node[concept] {(bio-) diversity measures}}
            child[grow=300 ] {node[concept] {volume etc.\ of subsets of $\R^n$}};    
        \end{tikzpicture}
    \end{center}
    \caption{Some of the rich connections of magnitude to various areas of mathematics}
    \label{fig:magnitude-connections}
\end{figure}
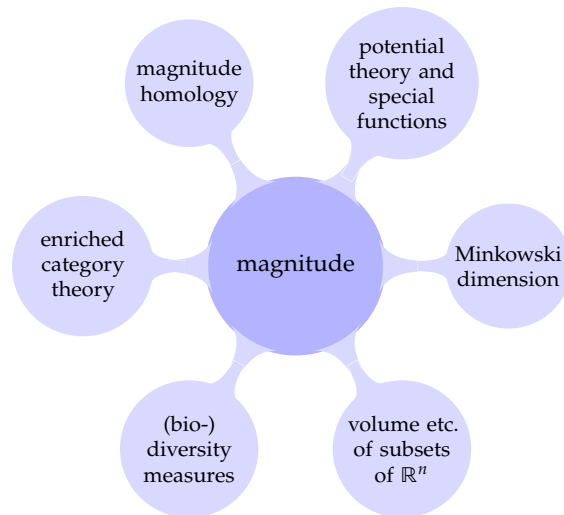

Magnitude has since been extensively studied and found to have connections to many areas of mathematics, some of these are illustrated in Figure~\ref{fig:magnitude-connections}.  (Leinster maintains a bibliography~\cite{Leinster:magnitude-bibliography} of magnitude-related papers; at the time of writing it has nearly~100 entries.)   Here I will just mention two of the related areas which are more connected to finite metric spaces, namely magnitude homology and further connections with bio-diversity.

\subsubsection*{Magnitude homology}
If you have a numerical invariant of something which is a natural number, (maybe a coefficient in a polynomial) then it sometimes arises as the dimension of a vector space and that vector space is an interesting and richer invariant than the numerical invariant, possibly having deeper structure such as functoriality.  A generalization of taking dimension of a vector space is taking the Euler characteristic of a sequence of vector spaces.  If you have a finite sequence of vector spaces $V_0, V_1, \dots, V_n$ then its Euler characteristic is the alternating sum of dimensions $\sum_i (-1)^i \dim V_i$.  Sometimes integer invariants are actually the Euler characteristic of some sequence of vector spaces.  The process of starting with an integer invariant and finding a richer invariant taking values in indexed sets of vector spaces is known as \define{categorification} -- the term is used in other contexts as well.  

A fundamental example is the following: the homology groups of a surface form a categorification of the classical Euler characteristic of the surface.  Another example comes from the Jones polynomial of knots and links; the Jones polynomial is a polynomial with integer coefficients; Khovanov~\cite{Khovanov:categorification} showed that each coefficient had a categorification to a sequence of vector spaces and these fitted together to form a ``bigraded homology theory'' which categorified the Jones polynomial.

Leinster~\cite{Leinster:graph-magnitude} showed that if you started with a graph (such as bipartite graph \(Y\) in Figure~\ref{subfig:bad-metric-space}) then by writing \(q = e^{-t}\) you could express the magnitude function of the graph as an integer power series in \(q\), so \(\sum a_j q^j\) for \(a_j\in \Z\).  Richard Hepworth and I~\cite{HepworthWillerton:magnitude-homology} then showed that each coefficient \(a_j\) could be categorified and these fitted together to form a bigraded homology theory of finite graphs, \define{magnitude homology}, which categorified the magnitude function.  Leinster and Shulman~\cite{LeinsterShulman:magnitude-homology} went on to show that there is a generalization to all metric spaces, which for finite metric spaces categorifies the magnitude function.

\subsubsection*{Diversity measures}

As mentioned above, Solow and Polasky originally discovered magnitude when devising a `utilitarian' diversity measure where a population is modelled by a metric space with the points corresponding to the species and the distances corresponding to the difference between the species.  Another standard approach is to model a population as a probability space where again the points correspond to species and the probabilities correspond to relative abundances.  One can then use entropy-like invariants of probability spaces to come up with diversity measures, in particular there is a one-parameter family of so-called Hill numbers, parametrized by \(q\in [0, \infty]\), with specific values of \(q\) giving (essentially) some well-known diversity measures such as Shannon entropy or Gini-Simpson diversity.  For every value of \(q\) the maximum value of the diversity is the number of species and this is obtained when the probability distribution is the uniform distribution, i.e., all species are equally abundant.  For details, a good mathematical reference is Leinster's book on entropy and diversity~\cite[Chapter~4]{Leinster:entropy-book} (the material below is covered in Chapter~6 of the book).

Tom Leinster and the mathematical biologist Christina Cobold~\cite{Leinster-Cobbold:species-similarity} showed that there was a way to combine these approaches by modelling a population as a metric space with a probability distribution on it and defining a one-parameter family of `similarity senstive diversity measures' which generalizes the classical Hill numbers.  The classical Hill numbers are recovered when the distance between each pair of species is taken to be \(\infty\).  Moreover, Leinster and Meckes~\cite{LeinsterMeckes:maximizing-diversity} showed that for a fixed metric space \(X\), every member of this family of diversity members has maximum value equal to a slight tweak of the magnitude of \(X\) and this is obtained when the probability distribution is a normalized version of a `weighting' (this is something used in the general definition of magnitude).

\section{Third Vignette: the Legendre-Fenchel transform}
Here we follow the same pattern as the previous two vignettes: the first subsection introduces the Legendre-Fenchel transform from a classical perspective and the second subsection brings in the enriched category perspective.

\subsection{The classical Legendre-Fenchel transform}

The Legendre-Fenchel transform occurs in various areas of mathematics, in particular in convex analysis and optimization where it is sometimes known as the conjugate transform (see, for example, Rockafellar's classic text~\cite[Section~12]{Rockafellar:convex-analysis}, or for a different approach, Touchette's unpublished notes~\cite{Touchette:Legendre-fenchel}).

\begin{figure}[tbh]
    \vskip -2em 
    \centering
    \(
        \begin{tikzpicture}[baseline=0.8cm]
            \begin{axis}[
                tick label style={font=\tiny},
                width = 10cm, height =3cm,
                axis equal image=true,
                axis x line=middle, axis y line = middle,
                xmin=-1.8,ymin=-0.3, 
                xmax=1.8, ymax=2, 
                xtick={-1,1}, xticklabels={$-1$,1},
                ytick={1}, yticklabels={1},
                x axis line style={style = -},
                y axis line style={style = -},
                xlabel=$x$, ylabel=\(z\), 
                yscale=1.5,xscale=1.5,
                legend style={at={(2,0.1)},anchor=south east}
                ]
                \addplot[mark=none, samples=60, brown,  thick,domain=-1.8:1.8] expression {x^4-2*x^2+1};
                \addplot[mark=none, samples=60, red,  thick,domain=-1.7:1.7] expression {-45/16*x-819/256};
            \end{axis}
        \end{tikzpicture}
        ~\xmapsto{\fenchb}~
        %
        %
        %
        %
        %
        %
        %
        %
        %
        %
        %
        \begin{tikzpicture}[baseline=0.8cm]
            \begin{axis}[
                    tick label style={font=\tiny},
                    width = 10cm, height =3cm,
                    axis equal image=true,
                    axis x line=middle, axis y line = middle,
                    xmin=-5.5 ,ymin=-2, 
                    xmax=5.5, ymax=7.5, 
                    xtick={-5,0,5}, xticklabels={-5,0, 5},
                    ytick={0, 5}, yticklabels={$0$, 5},
                    x axis line style={style = -},
                    y axis line style={style = -},
                    xlabel=$k$, ylabel=\(z\), 
                    yscale=1.5, 
                    xscale=1.5,
                ]
                \addplot[mark=none,brown, thick,domain=-5:5] coordinates {
                    (-5.0, 6.0820622169727505) (-4.583333333333333, 5.511553937108714) (-4.1666666666666661, 4.950627169826064) (-3.75, 4.3997009584123035) (-3.333333333333333, 3.859245500452226) (-2.9166666666666661, 3.329793158684577) (-2.5, 2.811953066896681) (-2.083333333333333, 2.3064309888580405) (-1.6666666666666661, 1.8140571285009757) (-1.25, 1.3358265013693784) (-0.83333333333333215, 0.872960229891512) (-0.41666666666666607, 0.4270041304350509) (0, 0) (0.41666666666666607, 0.4270041304350509) (0.83333333333333215, 0.872960229891512) (1.25, 1.3358265013693784) (1.6666666666666661, 1.8140571285009757) (2.083333333333333, 2.3064309888580405) (2.5, 2.811953066896681) (2.9166666666666661, 3.329793158684577) (3.333333333333333, 3.859245500452226) (3.75, 4.3997009584123035) (4.1666666666666661, 4.950627169826064) (4.583333333333333, 5.511553937108714) (5.0, 6.0820622169727505) 
                    };
                \addplot[mark=*, mark options={scale=0.4}, red, thick] coordinates {(-2.8125, 3.19921875)};
        \end{axis}
        \end{tikzpicture}
        ~\xmapsto{\rfenchb}~
        \begin{tikzpicture}[baseline=0.8cm]
            \begin{axis}[
                tick label style={font=\tiny},
                width = 10cm, height =3cm,
                axis equal image=true,
                axis x line=middle, axis y line = middle,
                xmin=-1.8, ymin=-0.3, 
                xmax=1.8, ymax=2, 
                xtick={-1,1}, xticklabels={$-1$,1},
                ytick={1}, yticklabels={1},
                x axis line style={style = -},
                y axis line style={style = -},
                xlabel=$x$, ylabel=\(z\), 
                yscale=1.5,
                xscale=1.5,
                legend style={at={(2,0.1)},anchor=south east}
                ]
                \addplot[mark=none, brown, thick,domain=-1:1] expression {0};
                \addplot[mark=none, brown, thick,domain=-1.8:-1] expression {x^4-2*x^2+1};
                \addplot[mark=none, brown, thick,domain=1:1.8] expression {x^4-2*x^2+1};
            \end{axis}
        \end{tikzpicture}
    \)
    \caption{Plots of the function \(f\colon x\mapsto (x^2-1)^2\), of \(\fenchbop{f}\) and of \(\rfenchb\circ \fenchbop{f}\).  A supporting hyperplane is pictured on the first plot, and a corresponding point is pictured on the second plot.}
    \label{fig:LF_example}
\end{figure}
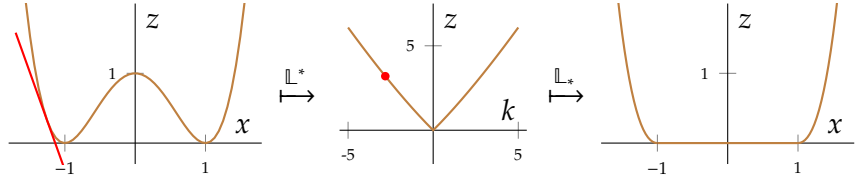

We fix a finite dimensional, real vector space \(V\) and we denote its linear dual by \(\dual{V}\).  The Legendre-Fenchel transform consists of a pair of functions between function spaces, where the functions in the function spaces take values in the extended real numbers, \([-\infty, +\infty]\coloneq \R\cup\{-\infty, +\infty\}\):
\[
\begin{tikzcd}[ampersand replacement=\&]
	{\Fun\bigl({V},[-\infty, +\infty]\bigr)} \& {\Fun\bigl(\dual{V},[-\infty, +\infty]\bigr).}
	\arrow["\fenchb", shift left, curve={height=-2pt}, from=1-1, to=1-2]
	\arrow["{\rfenchb }", shift left, curve={height=-2pt}, from=1-2, to=1-1]
\end{tikzcd}
\]
These functions are defined as follows:
\[
    \fenchbop{f}(k)\coloneqq \sup_{x\in V}\big\{\langle k,x\rangle -f(x)\big\},
    \quad 
    \rfenchbop{g}(x) \coloneqq \sup_{k\in \dual{V}}
    \big\{ \langle k, x \rangle -g(k) \big\}.
\]
You can see these illustrated in Figure~\ref{fig:LF_example} for the function \(f \colon \R \to [-\infty, +\infty]\), \(x\mapsto (x^2-1)^2\).

One way of interpreting the Legendre-Fenchel transform is via the notion of supporting hyperplanes.  Suppose we have a function \(f\colon V\to [-\infty, +\infty]\).  Consider the graph of \(f\) as a subset of \(V \times  [-\infty, +\infty]\), as is done in Figure~\ref{fig:LF_example}.  A \define{supporting hyperplane} \(H\) for \(f\) is a hyperplane in \(V \times  \R \subset V \times  [-\infty, +\infty]\) such that \(H\) touches and lies below the graph of \(f\).  A supporting hyperplane is illustrated in the first graph of Figure~\ref{fig:LF_example}.  If \(f\) is smooth then a supporting hyperplane is a tangent plane, but a tangent plane is not always a supporting hyperplane.  
For example, the tangent line to the graph at the point \((0, 1)\) in the first picture of Figure~\ref{fig:LF_example} is not a supporting hyperplane.  You can also observe that the graph in the second picture is not smooth at the origin, and there is a family of supporting hyperplanes which touch the graph there.

No supporting hyperplane can be vertical and so every supporting hyperplane \(H\) has a well-defined slope \(k_H\) which is a linear functional in \(\dual{V}\).  Moreover, there can be at most one supporting hyperplane with a given slope.  For a given a linear functional \(k \in \dual{V}\), if there is a supporting hyperplane \(H_k\) of \(f\) with slope \(k\) then the value of the Legendre-Fenchel transform \(\fenchbop{f}(k)\) is equal to minus the \(z\)-value of the intersection of the hyperplane \(H_k\) with the \(z\)-axis.  In Figure~\ref{fig:LF_example}, the pictured supporting hyperplane in the first picture is the tangent line to the graph at \(x=5/4\), this has equation \(z = -45/16x-819/256\); thus \(\fenchbop{f}(-45/16) = 819/256\), and this gives the point marked on the second graph.

On the other hand, if there is no supporting hyperplane of \(f\) with slope \(k\in \dual{V}\) then \(\fenchbop{f}(k)= +\infty\), \emph{unless} \(f\) is a specific function \(f_{+\infty}\) in which case \(\fenchbop{f_{+\infty}}(k)= -\infty\) for all \(k\); I will let you figure out which function \(f_{+\infty}\) is. 

This all tells us that the Legendre-Fenchel transform encodes the information about the supporting hyperplanes, so it is probably not surprising that you cannot always recover a function \(f\) from its Legendre-Fenchel transform, but that you can recover the envelope of the supporting hyperplanes, that is to say the closed convex hull of the function \(f\).  Here a closed convex function is one which is lower semi-continuous and such that if it takes the vales \(-\infty\) anywhere then it is the constant function at \(-\infty\).

What is perhaps surprising is that it is essentially applying the Legendre-Fenchel transform again, in the form of \(\rfenchb\) that recovers the closed convex hull of \(f\).  See the example in Figure~\ref{fig:LF_example}.  
This means that the composites \(\rfenchb\circ\fenchb\) and \(\fenchb\circ\rfenchb\) are \define{closed convex hull} operators.   In particular, \(\rfenchb\) and \(\fenchb\) become inverse isomorphisms when restricted to the sets of closed convex functions:
\[ 
    \Cvx\bigl({V},[-\infty, +\infty]\bigr)
    \cong 
    \Cvx\bigl(\dual{V}, [-\infty, +\infty]\bigr).
\]
This isomorphism is sometimes known as Fenchel-Moreau duality.

\subsection{The \texorpdfstring{\(\Rbar\)}{R}-category perspective}

(A general reference for the material here is my paper~\cite{Willerton:legendre-fenchel}.)
The key to understanding the Legendre-Fenchel transform in category theory terms is to work with the monoidal category \((\Rbar, +, 0)\).  Here \(\Rbar\) is like \(\Rplusbar\) but includes negative reals as well, so the object set is the extended real numbers, \(\ob\Rbar = [-\infty, +\infty]\), and there is a unique morphism \(a\to b\) precisely when \(a\ge b\).  The symmetric monoidal structure comes from extending addition from \((-\infty, +\infty)\) to \([-\infty, +\infty]\).  The only non-obvious sum needing to be defined is \((+\infty) + (-\infty)\) and it transpires that if we are to define an \(\Rbar\)-category \(\selfenr{\Rbar}\) we will need to take
\[
    (+\infty) + (-\infty) = +\infty = (-\infty) + (+\infty).
\]
We  then have the hom-objects given by subtraction
\[
    \selfenr{\Rbar}(a, b)\coloneq b-a
\]
where subtraction is extended in the obvious way to infinite values except the following cases:
\[
    (+\infty) - (+\infty) = -\infty = (-\infty) - (-\infty).
\]
This resolves a problem that is often encountered in this area, namely that you cannot consistently give meaning to the expression \(\infty - \infty\).  The enriched category perspective makes it clear that you should take care to write either \((+\infty) - (+\infty)\) or \((+\infty) + (-\infty)\) as these are not equal.

Our interest is in \(\Rbar\)-categories and these can be thought of as asymmetric metric spaces with possibly negative distances.  The axioms that are satisfied are the triangle inequality and that self-distance can be either \(0\) or \(-\infty\).  If \(x\in \ob X\) has \(X(x, x) = -\infty\) then \(x\) is an `infinite point' and the distance to other points must be either \(+\infty\) or \(-\infty\).

Now, getting to the context of the Legendre-Fenchel transform, if \(V\) is a finite-dimensional, real vector space then we will consider it as a discrete \(\Rbar\)-category meaning that the hom-objects, or distances, are given as follows:
\[
    V(x, x') \coloneq 
    \begin{cases}
        0 &\text{if\ } x = x'\\
        +\infty &\text{if\ } x \ne x'.
    \end{cases}
\]
We similarly consider the dual space \(\dual{V}\) as a discrete \(\Rbar\)-category.

The natural pairing between a vector space and its linear dual then gives rise to an \(\Rbar\)-profunctor:
\[
    \fenchpro\colon V \otimes \dual{V} \to \selfenr{\Rbar};
    \quad
    (x, k) \mapsto k(x).
\]
(As \(V\) is discrete we have \(V=V^\op\).)
From this we immediately get the Isbell-type adjunction~\eqref{eq:isbell-type-adjunction} coming from the profunctor \(\fenchpro\):
\[
\begin{tikzcd}[ampersand replacement=\&]
	{\funcat[\big]{V}{\selfenr{\Rbar}}} \& {\funcat[\big]{\dual{V}}{\selfenr{\Rbar}}^\op}
	\arrow["\fenchb", shift left, curve={height=-2pt}, from=1-1, to=1-2]
	\arrow["{\rfenchb }", shift left, curve={height=-2pt}, from=1-2, to=1-1]
\end{tikzcd}.
\]
On the level of a pair of functions between the sets of objects, this is precisely the Legendre-Fenchel transform; however, this comes with other structure built in.  For instance, both of \(\fenchb\) and \(\rfenchb\) are \(\Rbar\)-functors meaning that they are distance non-increasing with respect to the asymmetric, possibly negative distance on the function spaces given, for example, by
\[
    \funcat[\big]{V}{\selfenr{\Rbar}}(f,f')\coloneqq \sup_{x\in V} \{f'(x)-f(x)\}.
\]
The two \(\Rbar\)-functors also form an adjunction.  The nucleus of this adjunction, the fixed points, are precisely the closed convex functions and we get an \(\Rbar\)-\emph{isometry} between closed convex functions on \(V\) and closed convex functions on \(\dual{V}\):
\[
\begin{tikzcd}[ampersand replacement=\&, column sep=tiny]
	{\cvx[\big]{V}{\selfenr{\Rbar}}} \& \cong \& {\cvx[\big]{\dual{V}}{\selfenr{\Rbar}}^\op.}
	\arrow["\fenchb", shift left, curve={height=-2pt}, from=1-1, to=1-3]
	\arrow["{\rfenchb }", shift left, curve={height=-2pt}, from=1-3, to=1-1]
\end{tikzcd}
\]
This isometric structure is not hard to prove, once you are aware of it, and is known as Toland-Singer duality, although it is not part of the standard treatment of the Legendre-Fenchel transform.

\begin{figure}[ht]
    \vskip -2em
    \centering
    \(
        \begin{tikzpicture}[baseline=1.5cm]
            \begin{axis}[
                    tick label style={font=\tiny},
                    width = 10cm, height =4cm,
                    axis equal image=true,
                    axis x line=middle, axis y line = middle,
                    xmin=-1.5,ymin=-1.5, ymax=3.5, xmax=4.9,
                    xtick={-1,1,2,3,4}, xticklabels={-1,1,2,3,4},
                    ytick={-1,0,1,2,3}, yticklabels={-1,0,1,2,3},
                    x axis line style={style = -},y axis line style={style = -},
                    xlabel=$x$,
                    yscale=1.5, xscale=1.5,
                ]
                \addplot[mark=none,black,very thick,domain=-3:0] expression {-x};
                \node[black,font=\tiny] at (axis cs:2.5,3)  {$f$};
                \addplot[mark=none,black,very thick,domain=0:5] expression {x};
                \node[brown,font=\tiny] at (axis cs:3,0.5)  {$f'$};
                \addplot[mark=none,brown,very thick,domain=-3:2] expression {-x+1};
                \addplot[mark=none,brown,very thick,domain=2:5] expression {x-3};
            \end{axis}
        \end{tikzpicture}
        ~\xmapsto{\fenchb}~
        \begin{tikzpicture}[baseline=1.5cm]
            \begin{axis}
                [
                    tick label style={font=\tiny},
                    width = 7cm, height=4cm,
                    axis equal image=true,
                    axis x line=middle, axis y line = middle,
                    xmin=-2.3,ymin=-1.5, ymax=3.5, xmax=2.7,
                    xtick={-1,1}, xticklabels={$-1$,$1$},
                    ytick={-1,1,2,3},yticklabels={$-1$,$1$,$2$,$3$},
                    x axis line style={style = -},y axis line style={style = -},
                    xlabel=$k$,
                    yscale=1.5, xscale=1.5,
                ]
                \node[black,font=\tiny]
                        at (axis cs:0.6,0.4)  {$\fenchbop{f}$};
                \addplot[mark=*, mark options={scale=0.5}, black, very thick] coordinates {(-1,0) (1,0)};
                \node[brown,font=\tiny]
                        at (axis cs:0.1,2.5)  {$\fenchbop{f'}$};
                \addplot[mark=*, mark options={scale=0.5}, brown,very thick] coordinates {(-1,-1) (1,3)};
            \end{axis}
        \end{tikzpicture}
    \)
    \caption{}
    \label{fig:convex-isometry}
\end{figure}

Let's illustrate this \(\Rbar\)-isometry with an example.  In Figure~\ref{fig:convex-isometry} we can see the graphs of two convex functions \(f, f' \colon \R \to \closedinterval{-\infty}{+\infty}\) and the graphs of their transforms.  You should be able to verify quickly the fact that the distances are preserved: 
\begin{align*}
    \funcat[\big]{V}{\selfenr{\Rbar}}(f,f')
    &= 1 = 
    \funcat[\big]{\dual{V}}{\selfenr{\Rbar}}^\op\bigl(\fenchbop{f},\fenchbop{f'}\bigr);
    \\[0.5em]
    \funcat[\big]{V}{\selfenr{\Rbar}}(f', f)
    &= 3 = 
    \funcat[\big]{\dual{V}}{\selfenr{\Rbar}}^\op\bigl(\fenchbop{f'},\fenchbop{f}\bigr).
\end{align*}
    
This work shows that  \(\Rbar\)-category theory is a natural language in which to express this part of convex analysis; indeed a lot of Rockafellar's classic text~\cite{Rockafellar:convex-analysis} is expressible in this language.  It seems that this category theory perspective is opening up new vistas; see the work of Hanks et al~\cite{HanksEtAl:modeling-model-predictive-control} and that of Stein and Samuelson~\cite{SteinSamuleson:toward-a-compositional-framework}.

\printbibliography[
  heading=bibintoc,
  title={References}
]

\end{document}

%% file: pstex/Rplusbar_metric.pstex_t
\begin{picture}(0,0)%
\includegraphics{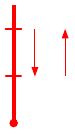}%
\end{picture}%
%
%
\setlength{\unitlength}{3947sp}%
\begingroup\makeatletter\ifx\SetFigFont\undefined%
\gdef\SetFigFont#1#2{%
  \fontsize{#1}{#2pt}%
  \selectfont}%
\fi\endgroup%
\begin{picture}(585,1057)(613,-331)
\put(940,262){\makebox(0,0)[lb]{\smash{{\SetFigFont{8}{9.6}{\color[rgb]{1,0,0}$0$}%
}}}}
\put(630,-283){\makebox(0,0)[rb]{\smash{{\SetFigFont{8}{9.6}{\color[rgb]{1,0,0}$0$}%
}}}}
\put(1183,263){\makebox(0,0)[lb]{\smash{{\SetFigFont{8}{9.6}{\color[rgb]{1,0,0}$b-a$}%
}}}}
\put(628, 77){\makebox(0,0)[rb]{\smash{{\SetFigFont{8}{9.6}{\color[rgb]{1,0,0}$a$}%
}}}}
\put(628,456){\makebox(0,0)[rb]{\smash{{\SetFigFont{8}{9.6}{\color[rgb]{1,0,0}$b$}%
}}}}
\end{picture}%

%% file: pstex/function_metric.pstex_t
\begin{picture}(0,0)%
\includegraphics{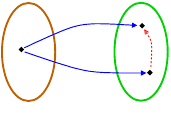}%
\end{picture}%
%
%
\setlength{\unitlength}{3947sp}%
\begingroup\makeatletter\ifx\SetFigFont\undefined%
\gdef\SetFigFont#1#2{%
  \fontsize{#1}{#2pt}%
  \selectfont}%
\fi\endgroup%
\begin{picture}(1356,992)(230,-396)
\put(466,-324){\makebox(0,0)[b]{\smash{{\SetFigFont{8}{9.6}{\color[rgb]{.75,.38,0}$X$}%
}}}}
\put(956,468){\makebox(0,0)[b]{\smash{{\SetFigFont{8}{9.6}{\color[rgb]{0,0,1}$g$}%
}}}}
\put(1359,-348){\makebox(0,0)[b]{\smash{{\SetFigFont{8}{9.6}{\color[rgb]{0,.69,0}$Y$}%
}}}}
\put(389,106){\makebox(0,0)[b]{\smash{{\SetFigFont{8}{9.6}{\color[rgb]{0,0,0}$x$}%
}}}}
\put(954, 84){\makebox(0,0)[b]{\smash{{\SetFigFont{8}{9.6}{\color[rgb]{0,0,1}$f$}%
}}}}
\end{picture}%

%% file: pstex/compact_subspaces.pstex_t
\begin{picture}(0,0)%
\includegraphics{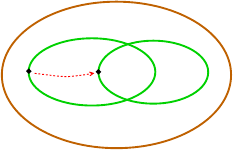}%
\end{picture}%
%
%
\setlength{\unitlength}{3947sp}%
\begingroup\makeatletter\ifx\SetFigFont\undefined%
\gdef\SetFigFont#1#2{%
  \fontsize{#1}{#2pt}%
  \selectfont}%
\fi\endgroup%
\begin{picture}(1864,1200)(169,-422)
\put(372,162){\makebox(0,0)[rb]{\smash{{\SetFigFont{8}{9.6}{\color[rgb]{0,0,0}$a$}%
}}}}
\put(1002,153){\makebox(0,0)[lb]{\smash{{\SetFigFont{8}{9.6}{\color[rgb]{0,0,0}$b$}%
}}}}
\put(732,-174){\makebox(0,0)[b]{\smash{{\SetFigFont{8}{9.6}{\color[rgb]{0,.69,0}$A$}%
}}}}
\put(1509,-172){\makebox(0,0)[b]{\smash{{\SetFigFont{8}{9.6}{\color[rgb]{0,.69,0}$B$}%
}}}}
\put(312,598){\makebox(0,0)[b]{\smash{{\SetFigFont{8}{9.6}{\color[rgb]{.75,.38,0}$M$}%
}}}}
\end{picture}%